\documentclass[11pt,oneside,reqno]{amsart}
\usepackage[T1]{fontenc}
\usepackage{lmodern}
\usepackage{amsmath,amssymb,amsthm,mathrsfs}
\usepackage{geometry}
\usepackage{xcolor}
\usepackage{tikz}
\usepackage[pagebackref=true,colorlinks=true,citecolor=blue,linkcolor=orange,urlcolor=blue]{hyperref}
\hypersetup{pdftitle={A spectral Moore problem and the Yang-Koolen conjecture},pdfauthor={Vishal Gupta}}
\newtheorem{theorem}{Theorem}[section]
\newtheorem{lemma}[theorem]{Lemma}
\newtheorem{proposition}[theorem]{Proposition}
\newtheorem{corollary}[theorem]{Corollary}

\theoremstyle{definition}

\newtheorem{question}[theorem]{Question}
\theoremstyle{remark}

\newcommand{\one}{\mathbf{1}}

\newcommand{\rev}[1]{{\hypersetup{linkcolor=blue,citecolor=blue,urlcolor=blue}\color{blue}#1}}
\DeclareMathOperator{\tr}{tr}
\begin{document}
\title[A spectral Moore problem]{A spectral Moore problem and the Yang-Koolen conjecture}
\author{Vishal Gupta}\thanks{Department of Mathematics, University of Rochester, Rochester, New York 14627, USA. Email: {vishalgupta@rochester.edu}.}

\begin{abstract}
We study the maximum order $v(k,\theta)$ of a connected $k$-regular graph whose second largest eigenvalue is at most $\theta$. For every real $\theta\geq1$, we prove that
\[
T(\theta) =\limsup_{k\to\infty}\bigl(v(k,\theta)-2k\bigr) =2\left\lfloor1+\frac{\theta^2}{4}\right\rfloor.
\]
In particular, $T(\theta)>2\theta$ for every real $\theta\geq4$, disproving the Yang-Koolen conjecture. We also determine $v(k,\theta)$ for all sufficiently large $k$ on explicit intervals of $\theta$. We prove a jump in the second largest eigenvalue above $1$. For each fixed $1\leq\theta<2$ and all sufficiently large $k$, every connected $k$-regular graph of order at least $2k+2$ whose second largest eigenvalue is at most $\theta$ is the crown graph. For $\theta=\sqrt2$, we prove that $v(k,\sqrt2)=2k+2$ for every $k\geq13$. Finally, we give a new upper bound on $v(k,\theta)$.
\end{abstract}

\maketitle

\section{Introduction}
Let $G$ be a connected $k$-regular graph on $n$ vertices and let $k=\lambda_1(G)>\lambda_2(G)\geq\cdots\geq\lambda_n(G)$ be its adjacency eigenvalues. The Alon-Boppana theorem \cite{AlonBoppana} states that,
for fixed $k\geq3$ and a sequence of connected $k$-regular
graphs $G_i$ with $|V(G_i)|\to\infty$,
\[
\liminf_{i\to\infty}\lambda_2(G_i)\geq2\sqrt{k-1}.
\]
Friedman \cite{Friedman1993} obtained quantitative lower bounds in terms of the diameter. Serre \cite{Serre} strengthened the Alon-Boppana theorem by showing that, for every $\varepsilon>0$, a positive proportion of the eigenvalues of a $k$-regular graph are at least $2\sqrt{k-1}-\varepsilon$, where the proportion depends only on $k$ and $\varepsilon$. Thus, for fixed $k\geq3$ and $\theta<2\sqrt{k-1}$, the orders of connected $k$-regular graphs with $\lambda_2\leq\theta$ are bounded. It is natural to ask for the maximum order of such a graph. This problem is known as a spectral Moore problem.

Let $v(k,\theta)$ denote the maximum order of a connected $k$-regular graph whose second largest eigenvalue is at most $\theta$. Richey, Shutty, and Stover \cite{RSS} studied this question using an effective form of Serre's argument and made several conjectures for particular values of $k$ and $\theta$. Marcus, Spielman, and Srivastava \cite{MSS} proved the existence of arbitrarily large bipartite Ramanujan graphs of every degree $k\geq3$. Their result implies that $v(k,\theta)$ is infinite when $\theta\geq2\sqrt{k-1}$. Therefore, $2\sqrt{k-1}$ is the precise boundary between the finite and infinite cases for each fixed $k\geq3$.

Cioab\u{a}, Koolen, Nozaki, and Vermette
\cite{Sebisiamvktheta} used the linear programming bound of Nozaki \cite{NozakiLP} to obtain a general upper bound on
$v(k,\theta)$. They determined several exact values and resolved a number of the conjectures in \cite{RSS}.
We recall their bound in Section~\ref{sec:prelim}.
The corresponding problem for bipartite regular graphs was studied by Cioab\u{a}, Koolen, and Nozaki \cite{CKN}. 
Cioab\u{a}, Koolen, Mimura, Nozaki, and Okuda \cite{CKMNO}
obtained bounds for bipartite biregular graphs through their work
on regular uniform hypergraphs. Lato \cite{Lato} later obtained
bounds that improve the earlier bounds in some cases.

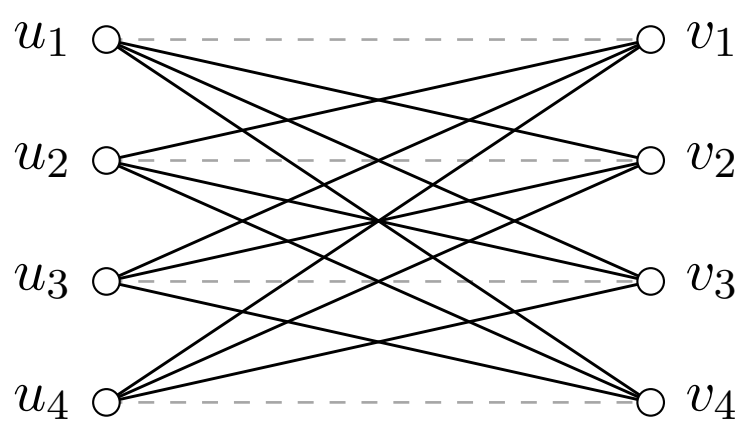
\begin{figure}[htbp]
\centering
\begin{tikzpicture}[x=1cm,y=1cm,
 vertex/.style={circle,draw=black,fill=white,inner sep=0pt,minimum size=5pt}]
\foreach \i in {1,...,4}{
  \coordinate (u\i) at (0,{-0.8*(\i-1)});
  \coordinate (v\i) at (3.6,{-0.8*(\i-1)});
}

\foreach \i in {1,...,4}{
  \draw[dashed,gray!70,line width=0.5pt] (u\i)--(v\i);
}
\foreach \i in {1,...,4}{
  \foreach \j in {1,...,4}{
    \ifnum\i=\j\relax\else
      \draw[black,line width=0.55pt] (u\i)--(v\j);
    \fi
  }
}
\foreach \i in {1,...,4}{
  \node[vertex,label=left:{$u_{\i}$}] at (u\i) {};
  \node[vertex,label=right:{$v_{\i}$}] at (v\i) {};
}
\end{tikzpicture}
\caption{The crown graph
$K_{4,4}\setminus\{u_iv_i:1\leq i\leq4\}$.
Solid lines are its edges. Dashed lines indicate the deleted
perfect matching.}
\label{fig:crown}
\end{figure}

We now fix $\theta\geq1$ and let $k$ tend to infinity.
For $k\geq2$, the crown graph is obtained by deleting a
perfect matching from $K_{k+1,k+1}$. It has valency $k$, order $2k+2$, and second largest eigenvalue $1$.
Hence, $v(k,\theta)\geq2k+2$.
Figure~\ref{fig:crown} shows the construction for $k=3$.

Yang and Koolen \cite{KoolenYang} proved that, for every
positive integer $\theta$, there is a constant $C(\theta)$
such that
\[
2k+2\leq v(k,\theta)\leq2k+C(\theta)
\]
for all $k>(\theta^2+4)/4$.
For a fixed real $\theta\geq1$, they defined
\begin{equation}\label{eq:Tdefinition}
T(\theta)=\limsup_{k\to\infty}\bigl(v(k,\theta)-2k\bigr)
\end{equation}
and conjectured that $T(\theta)=2\theta$ for every positive
integer $\theta$ \cite[Conjecture 6.3]{KoolenYang}.
For each fixed $\theta$, the condition
$\theta<2\sqrt{k-1}$ holds for all sufficiently large $k$,
so $v(k,\theta)$ is finite for these degrees.
Also, $v(k,\theta)\leq v(k,\lceil\theta\rceil)$.
Thus, their upper bound at $\lceil\theta\rceil$ shows that
the limsup in \eqref{eq:Tdefinition} is finite for every
real $\theta\geq1$.

In this paper, we determine $T(\theta)$ for every real
$\theta\geq1$. We prove that
\begin{equation}\label{eq:introT}
T(\theta)=2\left\lfloor1+\frac{\theta^2}{4}\right\rfloor.
\end{equation}
In particular, $T(\theta)>2\theta$ for every real
$\theta\geq4$. Consequently, the Yang-Koolen conjecture
holds for the integers $1,2,3$ and fails for every integer
$\theta\geq4$.
Moreover, for each fixed real $\theta\geq1$, we prove that
\[
v(k,\theta)\leq2k+2\left\lfloor1+\frac{\theta^2}{4}\right\rfloor
\]
for all sufficiently large $k$.
Equality holds for every sufficiently large $k$ whenever
$\theta$ is not of the form $2\sqrt{d-1}$ for an integer
$d\geq3$. At these remaining values of $\theta$, our
constructions attain the bound for infinitely many degrees. The lower bounds use bipartite complements of regular Ramanujan graphs. For the upper bound, we use Balla's result on equiangular lines \cite{Iballa} to obtain a graph of bounded maximum degree by Seidel switching. Regularity of the original graph gives a lower bound on the average degree of the switched graph. Jiang's bound on the spectral radius of a ball \cite{Jiang} gives an upper bound. Combining these bounds proves \eqref{eq:introT}.

We also study jumps in the possible values of $\lambda_2$ above $1$. Cioab\u{a}, Koolen, Nozaki, and Vermette \cite[Propositions 4.5 and 4.6]{Sebisiamvktheta} proved that
a connected cubic graph with $\lambda_2>1$ satisfies
$\lambda_2\geq\sqrt2$, while a connected $4$-regular graph
with $\lambda_2>1$ satisfies $\lambda_2\geq\sqrt5-1$.
The author, together with Cioab\u{a}, Nozaki, and Xiang
\cite[Proposition 4.3]{CGNX}, proved the same lower bound
$\sqrt5-1$ for connected $5$-regular graphs.
They also proved that every connected regular graph with
$\lambda_2>1$ satisfies
\[
\lambda_2>\sqrt\rho+\frac1{\sqrt\rho}-1\approx1.01980,
\]
where $\rho$ is the unique real root of $x^3-x-1=0$
\cite[Proposition 4.1]{CGNX}.

For graphs of order at least $2k+2$, we prove the
following jump result. For each fixed $1\leq\theta<2$
and all sufficiently large $k$, every connected $k$-regular
graph with $n\geq2k+2$ and $\lambda_2\leq\theta$ is the
crown graph. Thus, for each fixed $1<\theta<2$ and all sufficiently large $k$, no connected $k$-regular graph of order at least $2k+2$ has second largest eigenvalue in $(1,\theta]$. In particular, $v(k,\theta)=2k+2$ for each fixed $1\leq\theta<2$ and all sufficiently large $k$. We also give a direct consequence of Hoffman's theorem \cite{Hoffman} for $1\leq\theta<\sqrt2$.

Section~\ref{sec:bounds} gives a general upper bound on $v(k,\theta)$. Section~\ref{sec:lower} gives the bipartite construction, and Section~\ref{sec:exact} proves \eqref{eq:introT} and determines $v(k,\theta)$ when $k$ is sufficiently large. The jump results are proved in Section~\ref{sec:jump}. In Section~\ref{sec:sqrt2}, we use the classification of Gossett, Jiang, Teets, and Wellner \cite{GJTW} to prove that $v(k,\sqrt2)=2k+2$ for every $k\geq13$. We conclude with some open problems.

\section{Preliminaries}\label{sec:prelim}
We write $A_G$ for the adjacency matrix of a graph $G$, $J$ for the all-ones matrix, and $\one$ for the all-ones vector. The spectral radius of $G$ is denoted by $\rho(G)$, and its average degree by $\overline d(G)$. For a regular graph, its complement satisfies $\lambda_{\min}(\overline G)=-1-\lambda_2(G).$ This follows from $A_{\overline G}=J-I-A_G$,
whose eigenvalues on $\one^\perp$ are
$-1-\lambda_i(G)$ for $2\leq i\leq n$. 

We recall the bound of Cioab\u{a}, Koolen, Nozaki, and Vermette \cite[Theorem 2.3]{Sebisiamvktheta}.
\begin{theorem}\label{thm: vktheta bound}
Let $k,t\geq3$ be integers and let $c>0$ be real. Let $T(k,t,c)$ be the $t\times t$ tridiagonal matrix with lower diagonal $(1,\ldots,1,c)$, upper diagonal $(k,k-1,\ldots,k-1)$, and constant row sum $k$. If $\theta=\lambda_2(T(k,t,c))$, then
\[
v(k,\theta)\leq1+\sum_{i=0}^{t-3}k(k-1)^i+\frac{k(k-1)^{t-2}}{c}.
\]
Equality holds if and only if there is a distance-regular graph of valency $k$ and second largest eigenvalue $\theta$ whose distance-partition has quotient matrix $T(k,t,c)$.
\end{theorem}

We will use the following inequality from \cite[Theorem 1.15]{Iballa} to prove a new upper bound on $v(k, \theta)$.
\begin{theorem}\label{thm: Iballa bound}
Let $G$ be a $k$-regular graph on $n$ vertices with second largest eigenvalue $\lambda_2$. If $k-\lambda_2<\frac n2,$ then
\begin{equation}\label{eq: valency bound}
2\left(k-\frac{(k-\lambda_2)^2}{n}\right)
\leq\frac{\lambda_2(\lambda_2+1)(2\lambda_2+1)}{1-\frac2n(k-\lambda_2)}-\lambda_2(3\lambda_2+1).
\end{equation}
\end{theorem}

We also recall the theorem of Yang and Koolen \cite[Theorem 1.2]{KoolenYang}.
\begin{theorem}\label{thm: KoolenYang}
Let $\theta$ be a positive integer. Then there exists a constant $C(\theta)$ such that
\[2k+2\leq v(k,\theta)\leq2k+C(\theta)\]
holds for all $k>\frac{\theta^2+4}{4}$.
\end{theorem}
Note that $v(k, \theta)\leq v(k, \lceil \theta \rceil)$. Thus, $v(k,\theta)-2k$ is bounded by a constant depending on $\theta$ for all sufficiently large $k$.

\section{An upper bound on $v(k,\theta)$}\label{sec:bounds}

\begin{theorem}\label{thm:quadratic}
Let $\theta > 1$. If $k > \theta^3$, then
\begin{equation}\label{eq: quadratic upper bound}
    v(k, \theta) \leq \frac{(k-\theta) \left( 3k + 3\theta^2 + \sqrt{(3k + 3\theta^2)^2 - 8(k-\theta)(k-\theta^3)} \right)}{2(k - \theta^3)}.
\end{equation}
\end{theorem}
\begin{proof}
The condition $k>\theta^3$ implies $k\geq2$ and
$\theta<2\sqrt{k-1}$, so $v(k,\theta)$ is finite.
The classification in \cite[Theorem 3.2]{Sebisiamvktheta}
provides a graph attaining $v(k,1)$ with second largest
eigenvalue $1$ for every $k\geq2$.
If $v(k,\theta)>v(k,1)$, every graph attaining
$v(k,\theta)$ has second largest eigenvalue greater than $1$.
If $v(k,\theta)=v(k,1)$, we may choose one of these
extremal graphs with second largest eigenvalue $1$.

Thus, let $G$ be a connected $k$-regular graph of order
$n=v(k,\theta)$ with $1\leq\lambda_2\leq\theta.$ The crown graph gives $n\geq2k+2$.
Consequently, $0<k-\lambda_2<\frac n2.$ By Equation~\eqref{eq: valency bound}, we have
\[2\left( k - \frac{(k-\lambda_2)^2}{n}\right)\leq \frac{\lambda_2(\lambda_2+1)(2\lambda_2+1)}{1-\frac{2}{n}(k-\lambda_2)} - \lambda_2(3\lambda_2+1).\]
Let $x = \frac{2(k-\lambda_2)}{n}$. Since $k-\lambda_2 < \frac{n}{2}$ and $k > \lambda_2$, we have $0 < x < 1$. Substituting $x$ into the inequality and rearranging gives
\begin{equation}\label{eq: quadratic ineq}
   x^2(k-\lambda_2) - x(3k + 3\lambda_2^2) + 2(k - \lambda_2^3) \leq 0.
\end{equation}
Let $f(x) = ax^2 - bx + c$ represent the quadratic function in Equation (\ref{eq: quadratic ineq}). Since $k>\theta^3\geq \lambda_2^3\geq 1$, both $a$ and $c$ are positive. The discriminant
\begin{align*}
    \Delta &= (3k + 3\lambda_2^2)^2 - 8(k-\lambda_2)(k - \lambda_2^3)\\
    &= k^2 + (8\lambda_2^3 + 18\lambda_2^2+8\lambda_2)k + \lambda_2^4\\
    &> 0.
\end{align*}
Since $a,b$, and $c$ are all positive, by Vieta's formulas, we get that both the roots are positive. For the inequality $f(x) \leq 0$ to hold, $x$ must lie between the two roots. Therefore, $x$ is bounded below by the smaller root,
\[r_1 = \frac{3k + 3\lambda_2^2 - \sqrt{\Delta}}{2(k-\lambda_2)}.\]
Thus,
\[0< \frac{3k + 3\lambda_2^2 - \sqrt{\Delta}}{2(k-\lambda_2)} \leq x <1.\]
Now, substituting $x= \frac{2(k-\lambda_2)}{n}$ and solving for $n$ gives
\[ n \leq \frac{(k-\lambda_2) \left( 3k + 3\lambda_2^2 + \sqrt{\Delta} \right)}{2(k - \lambda_2^3)}:= F(\lambda_2). \] Finally, note that $F(\lambda_2)$ is an increasing function in $\lambda_2$ when $k>\lambda_2^3$. Therefore,
\[ v(k,\theta) = n \leq F(\theta) = \frac{(k-\theta) \left( 3k + 3\theta^2 + \sqrt{(3k + 3\theta^2)^2 - 8(k-\theta)(k-\theta^3)} \right)}{2(k - \theta^3)}. \]
\end{proof}

For fixed $\theta>1$, the right-hand side of \eqref{eq: quadratic upper bound} has the following asymptotic expansion. Let $a_\theta=8\theta^3+18\theta^2+8\theta$. Then
\[
\sqrt{k^2+a_\theta k+\theta^4}
=k+\frac{a_\theta}{2}+O_\theta(k^{-1}),
\]
and hence
\begin{align}
\frac{(k-\theta)(3k+3\theta^2+\sqrt{k^2+a_\theta k+\theta^4})}{2(k-\theta^3)}
&=\frac{4k^2+(4\theta^3+12\theta^2)k+O_\theta(1)}{2k-2\theta^3}\notag\\
&=2k+4\theta^3+6\theta^2+O_\theta(k^{-1}).\label{eq:asymptotic}
\end{align}
Thus, for fixed $\theta$, this bound has leading term $2k$.

For comparison, when $k\geq3$, $k>\theta^2$, and $\theta>1$, take $t=3$ and $c=(k-\theta^2)/(\theta+1)$ in Theorem~\ref{thm: vktheta bound}. The two eigenvalues other than $k$ of $T(k,3,c)$ are the roots of $x^2+cx+c-k$, and the larger is $\theta$. Therefore,
\begin{align*}
v(k,\theta)&\leq1+k+\frac{k(k-1)(\theta+1)}{k-\theta^2}\\
&=(\theta+2)k+\theta^3+\theta^2-\theta+O_\theta(k^{-1}).
\end{align*}
In particular, at $\theta=\sqrt2$, this gives
\[
v(k,\sqrt2)\leq\frac{(2+\sqrt2)k(k-1)-2}{k-2},
\]
whereas \eqref{eq:asymptotic} gives $2k+8\sqrt2+12+O(k^{-1})$. We determine $T(\theta)$ in Section~\ref{sec:exact}.

\section{Bipartite complements and a lower bound on $T(\theta)$}\label{sec:lower}

\begin{lemma}\label{lem:complement}
Let $H$ be an $s$-regular bipartite graph with parts of size $m\geq2$, and let $G=K_{m,m}\setminus E(H)$. Then $G$ is $(m-s)$-regular. Apart from the two trivial eigenvalues, the adjacency spectra of $H$ and $G$ agree as multisets. In particular, if $m>2s$, then $G$ is connected and $\lambda_2(G)=\lambda_2(H)$.
\end{lemma}
\begin{proof}
Let
\[
A_H=\begin{pmatrix}0&B\\B^T&0\end{pmatrix}
\]
be the adjacency matrix of $H$, where $B$ is its biadjacency matrix. Then
\[
A_G=\begin{pmatrix}0&J-B\\J-B^T&0\end{pmatrix}
\]
is the adjacency matrix of $G$. $G$ is $(m-s)$-regular. The vectors $(\one,\one)$ and $(\one,-\one)$ give the two trivial eigenvalues $s$ and $-s$ for $H$, and $m-s$ and $-(m-s)$ for $G$. The remaining eigenvectors belong to the orthogonal complement
\[
W=\{(x,y)\in\mathbb R^{2m}\mid \one^Tx=0,\ \one^Ty=0\}.
\]
Let $(\lambda,(x,y))$ be an eigenpair of $H$ with $(x,y)\in W$. Then
\[
By=\lambda x,\qquad B^Tx=\lambda y.
\]
Since $Jx=Jy=0$, we have
\[
A_G\begin{pmatrix}x\\-y\end{pmatrix}
=\begin{pmatrix}-(J-B)y\\(J-B^T)x\end{pmatrix}
=\begin{pmatrix}By\\-B^Tx\end{pmatrix}
=\lambda\begin{pmatrix}x\\-y\end{pmatrix}.
\]
Thus, $(\lambda,(x,-y))$ is an eigenpair of $G$. This proves that the remaining eigenvalues agree, including multiplicities. Finally, let $k=m-s$. Each component of a bipartite $k$-regular graph has at least $k$ vertices in each part. If $m>2s$, then $2k>m$, so two components are impossible. Thus, $G$ is connected.
\end{proof}

\begin{theorem}\label{thm:lower}
For every integer $s\geq3$ and every real $\theta\geq2\sqrt{s-1}$, there are infinitely many degrees $k$ for which
\[v(k,\theta)\geq2k+2s.\]
Consequently, for every real $\theta\geq1$,
\begin{equation}\label{eq:lowerT}
T(\theta)\geq2\left\lfloor1+\frac{\theta^2}{4}\right\rfloor.
\end{equation}
\end{theorem}
\begin{proof}
By Marcus, Spielman, and Srivastava \cite{MSS}, there are connected $s$-regular bipartite Ramanujan graphs $H_i$ whose two parts have size $m_i$, where $m_i\to\infty$. Let $G_i$ be the bipartite complement of $H_i$. For all sufficiently large $i$, we have $m_i>2s$, so Lemma~\ref{lem:complement} shows that $G_i$ is connected and has degree $k_i=m_i-s$. Its order is
\[
|V(G_i)|=2m_i=2k_i+2s,
\]
and
\[
\lambda_2(G_i)=\lambda_2(H_i)\leq2\sqrt{s-1}\leq\theta.
\]
Therefore $v(k_i,\theta)\geq2k_i+2s$. Since $k_i\to\infty$, this proves the first assertion and gives $T(\theta)\geq2s$.

Now fix $\theta\geq1$ and let $L=\lfloor1+\theta^2/4\rfloor$. We prove $T(\theta)\geq2L$ by considering the possible values of $L$. If $L\geq3$, take $s=L$. From $L\leq1+\theta^2/4$, we get $2\sqrt{L-1}\leq\theta.$ Thus, $T(\theta)\geq2L$. If $L=2$, then $\theta\geq2$. For every integer $m>4$, take $H=C_{2m}$. This graph is $2$-regular, has parts of size $m$, and has second largest eigenvalue
\[
\lambda_2(H)=2\cos\frac{\pi}{m}<2\leq\theta.
\]
By Lemma~\ref{lem:complement}, its bipartite complement $G$ is connected, has degree $k=m-2$, and has the same second largest eigenvalue. Since $|V(G)|=2m=2k+4,$ we get $v(k,\theta)-2k\geq4$ for every $k\geq3$. Hence $T(\theta)\geq4=2L$. Finally, if $L=1$, the crown graph has degree $k$, order $2k+2$, and second largest eigenvalue $1\leq\theta$ for every $k\geq2$. Thus $v(k,\theta)-2k\geq2$, and $T(\theta)\geq2=2L$. In each case, $T(\theta)\geq2L$, which proves \eqref{eq:lowerT}.
\end{proof}

\begin{corollary}\label{cor:disproof}
For every real $\theta\geq4$, $T(\theta)>2\theta$. In particular, the Yang-Koolen conjecture is false for every integer $\theta\geq4$.
\end{corollary}
\begin{proof}
Take $s=\lfloor\theta\rfloor+1$. Then $s>\theta$ and
\[
2\sqrt{s-1}=2\sqrt{\lfloor\theta\rfloor}\leq2\sqrt\theta\leq\theta.
\]
Theorem~\ref{thm:lower} gives $T(\theta)\geq2s>2\theta$.
\end{proof}

\section{The exact value of $T(\theta)$}\label{sec:exact}
Given a partition $V(G)=X\cup Y$, Seidel switching keeps the edges within $X$ and within $Y$, and interchanges edges and nonedges between $X$ and $Y$. The Seidel matrix of a graph $G$ is $S_G=J-I-2A_G$. If $D$ is diagonal with entries $1$ on $X$ and $-1$ on $Y$, the switched graph $F$ has Seidel matrix $S_F=DS_GD$. 

Let $\theta\geq1$ and let $\alpha=(2\theta+1)^{-1}$. Let $G$ be a $k$-regular graph on $n\geq2k+2$ vertices with $\lambda_2(G)\leq\theta$. Let
\begin{equation}\label{eq:Gram}
M=I+\alpha S_G
=(1-\alpha)I+\alpha J-2\alpha A_G
=2\alpha\left(\theta I+\frac12J-A_G\right).
\end{equation}
Choose an orthonormal eigenbasis $u_1,\ldots,u_n$ of $A_G$ with $u_1=\one/\sqrt n$. For $i\geq2$, we have $\one^Tu_i=0$, so
\[
Mu_i=2\alpha(\theta-\lambda_i(G))u_i,
\]
and
\[
M\one=2\alpha\left(\theta+\frac n2-k\right)\one.
\]
For $i\geq2$, we have $\lambda_i(G)\leq\lambda_2(G)\leq\theta$,
and $n\geq2k+2$ gives $\theta+n/2-k\geq\theta+1>0$.
Thus, all eigenvalues of $M$ are nonnegative, and $M$ is positive semidefinite. Hence, there are vectors $p_1,\ldots,p_n$ whose Gram matrix
is $M$. The diagonal entries of $M$ are $1$, so each $p_i$
is a unit vector. Its off-diagonal entries are $\alpha$
or $-\alpha$, so
\[
|\langle p_i,p_j\rangle|=\alpha
\qquad(i\neq j).
\]
Thus, the lines spanned by these vectors are equiangular
with common angle $\arccos\alpha$. Replacing selected vectors $p_i$ by $-p_i$ changes the Gram
matrix to $DMD$, where $D$ is diagonal with entry $-1$ at each selected index and $1$ elsewhere. Since $D^2=I$,
\[
DMD=I+\alpha DS_GD,
\]
this change corresponds to Seidel switching between the two sign classes.

\begin{lemma}\label{lem:switch}
Let $\theta\geq1$, and let $G$ be a $k$-regular graph on $n$ vertices with $n\geq2k+2$ and $\lambda_2(G)\leq\theta$. If
\[
n>\binom{4\theta(\theta+1)}2,
\]
then $G$ is switching equivalent to a graph $F$ satisfying
\[
\Delta(F)\leq D_\theta
:=\left\lceil
\theta^3+\frac32\theta^2+\frac32\theta
\right\rceil-1.
\]
Furthermore, every vector $x\in\mathbb R^n$ whose coordinates sum to zero satisfies
\begin{equation}\label{eq:compression}
x^TA_Fx\leq\theta\|x\|^2.
\end{equation}
\end{lemma}
\begin{proof}
Let $\alpha=(2\theta+1)^{-1}$, let $M=I+\alpha S_G$, and let 
$a=\alpha^{-2}-1=4\theta(\theta+1)$. Since $M\succeq0$, we have
\[
\tr(M^2)\leq\lambda_1(M)\tr(M).
\]
The entries of $M$ give $\tr(M)=n$ and
$\tr(M^2)=n(1+\alpha^2(n-1))$. Since $\alpha^2=1/(a+1)$ and $n>a(a-1)/2$, we obtain
\[
\lambda_1(M)\geq1+\alpha^2(n-1)
=\frac{a+n}{a+1}
>\frac{a+a(a-1)/2}{a+1}
=\frac a2
=\frac{1-\alpha^2}{2\alpha^2}.
\]

Let $y$ be an eigenvector of $M$ corresponding to $\lambda_1(M)$. Choose the diagonal sign matrix $D$ so that $Dy$ is nonnegative, taking a diagonal entry $-1$ precisely when the corresponding coordinate of $y$ is negative. Then
\[
(DMD)(Dy)=DMy=\lambda_1(M)Dy.
\]
Since $DMD$ and $M$ have the same eigenvalues, $Dy$ is a nonnegative eigenvector for the largest eigenvalue of $DMD$. Let $F$ be the corresponding switched graph. Using \cite[Lemmas 3.1 and 3.6]{Iballa} we get that
\[
\Delta(F)<\frac{1+3\alpha^2-4\alpha^3}{8\alpha^3}
=\theta^3+\frac32\theta^2+\frac32\theta.
\]
Since $\Delta(F)$ is an integer, the degree bound follows. Finally, let $x\in\mathbb R^n$ have coordinate sum zero. Then $Jx=0$. Since
\[
DMD=(1-\alpha)I+\alpha J-2\alpha A_F
\]
is positive semidefinite, we have
\[
0\leq x^TDMDx
=(1-\alpha)\|x\|^2-2\alpha x^TA_Fx.
\]
Rearranging and using $(1-\alpha)/(2\alpha)=\theta$ gives \eqref{eq:compression}.
\end{proof}

Regularity of $G$ gives the following relations for the degrees of the switched graph.

\begin{lemma}\label{lem:cut}
Let $G$ be a $k$-regular graph of order $n=2k+2s$, where $s\geq1$. Let $d_0<k$ be a number such that switching a partition $V(G)=X\cup Y$ gives a graph $F$ with $\Delta(F)\leq d_0$.
Let $t=(|X|-|Y|)/2$, so that
\[
|X|=\frac n2+t,\qquad |Y|=\frac n2-t,
\]
and let $e_X,e_Y$ be the numbers of edges in $G[X],G[Y]$, respectively. Then
\[
s+|t|\leq d_0,\qquad e_X-e_Y=tk,
\]
and
\begin{equation}\label{eq:averagecut}
\overline d(F)
=s-\frac{2t^2}{n}+\frac{4(e_X+e_Y)}n
\geq s+2|t|-\frac{2d_0^2}{n}
\geq\lceil s\rceil-\frac{2d_0^2}{n}.
\end{equation}
The numbers $s$ and $t$ are integers when $n$ is even
and half-integers when $n$ is odd.
\end{lemma}

\begin{proof}
Both parts are nonempty, since otherwise $F=G$,
contrary to $\Delta(F)\leq d_0<k$. For each vertex $v$, let $a_v$ be its number of neighbors in $G$ within its own part. Switching keeps these neighbors and interchanges neighbors and nonneighbors in the other part. Thus,
\[
d_F(v)=
\begin{cases}
|Y|-(k-a_v)+a_v=s-t+2a_v,&v\in X,\\
|X|-(k-a_v)+a_v=s+t+2a_v,&v\in Y.
\end{cases}
\]
Since $a_v\geq0$ and both parts are nonempty, $s-t\leq d_0$ and $s+t\leq d_0$. Hence, $s+|t|\leq d_0$. Counting the edges of $G$ between $X$ and $Y$ from the two parts gives
\[
k|X|-2e_X=k|Y|-2e_Y.
\]
Since $|X|-|Y|=2t$, this gives $e_X-e_Y=tk$. Summing the degree formulas, we obtain
\[
\begin{aligned}
n\overline d(F)
&=(s-t)|X|+(s+t)|Y|+4e_X+4e_Y\\
&=ns-2t^2+4(e_X+e_Y).
\end{aligned}
\]
This proves the equality in \eqref{eq:averagecut}.
Also, $e_X+e_Y\geq|e_X-e_Y|=|t|k.$ Using $k=n/2-s$, we therefore obtain
\[
\overline d(F)
\geq s+2|t|-\frac{2t^2+4s|t|}{n}.
\]
Since $s+|t|\leq d_0$,
\[
2t^2+4s|t|
=2(s+|t|)^2-2s^2
\leq2d_0^2.
\]
This gives the first inequality in \eqref{eq:averagecut}. If $n$ is even, then $s$ is an integer and $s+2|t|\geq\lceil s\rceil$. If $n$ is odd, then $s,t$ are half-integers and $|t|\geq1/2$, so again $s+2|t|\geq\lceil s\rceil$. This proves the final inequality.
\end{proof}

\begin{lemma}\label{lem:averageupper}
Let $F_i$ be graphs with orders $n_i\to\infty$ and uniformly
bounded maximum degree. Suppose that \eqref{eq:compression}
holds for every $F_i$, with a fixed $\theta\geq1$. Then
\[
\limsup_{i\to\infty}\overline d(F_i)\leq1+\frac{\theta^2}{4}.
\]
Equivalently, for every $\varepsilon>0$,
\[
\overline d(F_i)\leq1+\frac{\theta^2}{4}+\varepsilon
\]
for all sufficiently large $i$.
\end{lemma}
\begin{proof}
Let $d_0\geq1$ be a common maximum-degree bound, and let
$d_i=\overline d(F_i)$. Fix an integer $r\geq1$. If $d_i\geq1$, Lemma 6 of \cite{Jiang} supplies a ball $H_i$ of radius $r$ with
\[
\rho(H_i)\geq2\sqrt{d_i-1}\cos\frac{\pi}{r+2}.
\]
The order of $H_i$ is at most $B_r=1+d_0+\cdots+d_0^r$.
Let $x$ be a Perron vector of $H_i$ with $\|x\|=1$,
extended by zero to the remaining vertices of $F_i$.
Let
\[
c=\langle x,\one\rangle,\qquad
z=x-\frac{c}{n_i}\one.
\]
Then $\langle z,\one\rangle=0$,
$\|z\|^2=1-c^2/n_i\leq1$, and Cauchy-Schwarz gives
$c^2\leq|V(H_i)|\leq B_r$.
Since $x$ is nonnegative and
$A_{F_i}\one\leq d_0\one$ coordinatewise,
\[
\begin{aligned}
z^TA_{F_i}z
&=x^TA_{F_i}x-\frac{2c}{n_i}\one^TA_{F_i}x
  +\frac{c^2}{n_i^2}\one^TA_{F_i}\one\\
&\geq\rho(H_i)-\frac{2d_0c^2}{n_i}\\
&\geq\rho(H_i)-\frac{2d_0B_r}{n_i}.
\end{aligned}
\]
This is also the projected-vector estimate in
\cite[Lemma 3.9]{Iballa}. By \eqref{eq:compression},
$z^TA_{F_i}z\leq\theta\|z\|^2\leq\theta$. Therefore,
\[
2\sqrt{d_i-1}\cos\frac{\pi}{r+2}
\leq\theta+\frac{2d_0B_r}{n_i},
\]
and hence
\[
d_i\leq
1+\frac{(\theta+2d_0B_r/n_i)^2}
        {4\cos^2(\pi/(r+2))}.
\]

Let $\varepsilon>0$. Since $\cos(\pi/(r+2))\to1$ as
$r\to\infty$, we can first choose $r$ such that
\[
1+\frac{\theta^2}{4\cos^2(\pi/(r+2))}
<1+\frac{\theta^2}{4}+\varepsilon.
\]
For this fixed $r$, the preceding upper bound on $d_i$
tends to the left-hand side as $i\to\infty$.
Consequently,
\[
d_i\leq1+\frac{\theta^2}{4}+\varepsilon
\]
for all sufficiently large $i$ with $d_i\geq1$.
The indices with $d_i<1$ already satisfy this inequality.
Since $\varepsilon>0$ was arbitrary, the result follows.
\end{proof}

\begin{theorem}\label{thm:exactT}
For every fixed real $\theta\geq1$, all sufficiently large degrees $k$ satisfy
\[
v(k,\theta)\leq2k+2\left\lfloor1+\frac{\theta^2}{4}\right\rfloor.
\]
Moreover,
\[
T(\theta)=2\left\lfloor1+\frac{\theta^2}{4}\right\rfloor.
\]
\end{theorem}
\begin{proof}
Let
\[
C=1+\frac{\theta^2}{4},\qquad
L=\lfloor C\rfloor,\qquad
\varepsilon=\frac{L+1-C}{3}>0.
\]
Suppose that the upper bound fails for arbitrarily large
degrees. Then there are degrees $k_i\to\infty$ and connected
$k_i$-regular graphs $G_i$ with $\lambda_2(G_i)\leq\theta$
and orders $n_i>2k_i+2L$.

Since $n_i\to\infty$, Lemma~\ref{lem:switch} applies for
all sufficiently large $i$. Let $F_i$ be the resulting
switched graphs, and let $d_0=D_\theta$ be their common
maximum-degree bound. We may assume that $k_i>d_0$.
Let
\[
s_i=\frac{n_i-2k_i}{2}.
\]
Since $s_i>L$, we have $\lceil s_i\rceil\geq L+1$.
Lemma~\ref{lem:cut} therefore gives
\[
\overline d(F_i)
\geq\lceil s_i\rceil-\frac{2d_0^2}{n_i}
\geq L+1-\frac{2d_0^2}{n_i}.
\]
For all sufficiently large $i$, we have
$2d_0^2/n_i<\varepsilon$, and hence
\[
\overline d(F_i)>L+1-\varepsilon=C+2\varepsilon.
\]
On the other hand, Lemma~\ref{lem:averageupper} gives $\overline d(F_i)\leq C+\varepsilon$ for all sufficiently large $i$, a contradiction. Thus, the eventual upper bound holds, and $T(\theta)\leq2L$.
If $L=1$, crown graphs give the matching lower bound.
If $L=2$, the bipartite complements of $C_{2(k+2)}$ give
the matching lower bound. If $L\geq3$,
Theorem~\ref{thm:lower} gives the matching lower bound
for infinitely many degrees. Thus, in every case,
$T(\theta)\geq2L$. This proves the equality.
\end{proof}

\begin{corollary}\label{cor:plateaux}
For each fixed $1\leq\theta<2$, $v(k,\theta)=2k+2$ for all sufficiently large $k$.
For each fixed $2\leq\theta<2\sqrt2$, $v(k,\theta)=2k+4$ for all sufficiently large $k$.
\end{corollary}
\begin{proof}
Theorem~\ref{thm:exactT} gives the upper bounds.
Crown graphs attain the first lower bound.
For the second, take the bipartite complement of $C_{2(k+2)}$.
For $k>2$, it is connected, has degree $k$, order $2k+4$,
and second eigenvalue $2\cos(\pi/(k+2))<2$.
\end{proof}

\begin{corollary}
For a positive integer $q$, $T(q)=2q$ holds if and only if $q\in\{1,2,3\}$.
\end{corollary}
\begin{proof}
Theorem~\ref{thm:exactT} gives $T(1)=2$, $T(2)=4$, and $T(3)=6$.
For every integer $q\geq4$, we have $q^2/4\geq q$, so
\[
\left\lfloor1+\frac{q^2}{4}\right\rfloor\geq q+1.
\]
Hence, $T(q)\geq2q+2>2q$.
\end{proof}

We also determine $v(k,\theta)$ for all sufficiently large $k$.
\begin{corollary}\label{cor:noncritical}
Let $\theta\geq1$ be fixed, and let
$L=\lfloor1+\theta^2/4\rfloor$.
If $L\geq3$ and $\theta>2\sqrt{L-1}$, then $v(k,\theta)=2k+2L$ for every sufficiently large $k$.
\end{corollary}
\begin{proof}
Let
\[
\varepsilon=\frac12\min\{\theta-2\sqrt{L-1}, L-2\sqrt{L-1}\}>0.
\]
By Brito, Dumitriu, and Harris \cite[Theorem 3.2]{BDH},
for every sufficiently large $m$ there is a simple
$L$-regular bipartite graph $H$ with parts of size $m$
such that
\[
\lambda_2(H)\leq2\sqrt{L-1}+\varepsilon<\min\{\theta,L\}.
\]
The graph $H$ is connected, since a disconnected
$L$-regular graph has second largest eigenvalue $L$.
Let $m=k+L$, and take the bipartite complement of $H$.
For all sufficiently large $k$, we have $m>2L$, so
Lemma~\ref{lem:complement} gives a connected graph of degree
$k$, order $2k+2L$, and second largest eigenvalue less
than $\theta$. This proves $v(k,\theta)\geq2k+2L$
for every sufficiently large $k$.
Theorem~\ref{thm:exactT} gives the reverse inequality.
\end{proof}
\begin{corollary}\label{cor:denselower}
Let $d\geq2$ be a fixed integer, and let $G_i$ be connected $k_i$-regular graphs with $k_i\to\infty$ and $|V(G_i)|\geq2k_i+2d-1$. Then
\[
\liminf_{i\to\infty}\lambda_2(G_i)\geq2\sqrt{d-1}.
\]
The bound is attained in the limit by a sequence of graphs of order $2k_i+2d$.
\end{corollary}
\begin{proof}
For each $1\leq\eta<2\sqrt{d-1}$, we have $\lfloor1+\eta^2/4\rfloor\leq d-1$. Theorem~\ref{thm:exactT} therefore excludes $\lambda_2(G_i)\leq\eta$ for all sufficiently large $i$. Letting $\eta$ increase to $2\sqrt{d-1}$ proves the bound. Bipartite complements of the Ramanujan graphs in Theorem~\ref{thm:lower} give sharpness for $d\geq3$. For $d=2$, use the bipartite complements of even cycles.
\end{proof} 
\begin{corollary}\label{cor:least-order}
For every fixed real $\theta\geq1$, there exists an integer
$k_0(\theta)$ such that every $k$-regular graph $G$ of order $n$
satisfying $n-k-1\geq k_0(\theta)$ and $\lambda_{\min}(G)\geq-1-\theta$ has
\[
n\geq2k-2\left\lfloor\frac{\theta^2}{4}\right\rfloor.
\]
The bound is attained for infinitely many degrees $k$.
\end{corollary}

\begin{proof}
Choose $k_0(\theta)>\theta$ so that Theorem~\ref{thm:exactT}
applies to every degree at least $k_0(\theta)$.
Let $H=\overline G$ and $q=n-k-1$. Then $H$ is $q$-regular and
\[
\lambda_2(H)=-1-\lambda_{\min}(G)\leq\theta.
\]
Since $q>\theta$, the graph $H$ is connected.
Theorem~\ref{thm:exactT} therefore gives
\[
n\leq2q+2\left\lfloor1+\frac{\theta^2}{4}\right\rfloor
 =2(n-k-1)+2\left\lfloor1+\frac{\theta^2}{4}\right\rfloor.
\]
Rearranging proves the bound.

For sharpness, let $L=\lfloor1+\theta^2/4\rfloor$.
By Theorem~\ref{thm:exactT}, for infinitely many degrees $q$
there is a connected $q$-regular graph $H$ of order $2q+2L$
with $\lambda_2(H)\leq\theta$. Its complement has degree
$k=q+2L-1$, smallest eigenvalue at least $-1-\theta$, and order
\[
2q+2L=2k+2-2L
 =2k-2\left\lfloor\frac{\theta^2}{4}\right\rfloor.
\]
Here $n-k-1=q\to\infty$, so the hypothesis also holds for all
sufficiently large degrees in this sequence.
\end{proof}

\section{Jumps in the second largest eigenvalue}\label{sec:jump}
We first give a direct consequence of Hoffman's theorem \cite{Hoffman}. We use the following uniform formulation, stated in \cite[Theorem 1.1]{KYY}. For every $\alpha\in(-1-\sqrt2,-2]$, there exists an integer $D(\alpha)$ such that a connected graph $H$ with $\delta(H)\geq D(\alpha)$ and $\lambda_{\min}(H)\geq\alpha$ has $\lambda_{\min}(H)\geq-2$. The theorem in \cite{KYY} treats $\alpha\leq\lambda_{\min}(H)\leq-2$; the case $\lambda_{\min}(H)>-2$ already has the required conclusion.

\begin{proposition}\label{prop:hoffmanjump}
Fix $1\leq\theta<\sqrt2$. There exists an integer $d_\theta$ such that every regular graph $G$ of degree $k$ and order $n$ satisfying $n-k-1\geq d_\theta$ and $\lambda_2(G)\leq\theta$ has $\lambda_2(G)\leq1$.
\end{proposition}
\begin{proof}
Let $\alpha=-1-\theta\in(-1-\sqrt2,-2]$, and choose $d_\theta=D(\alpha)$. The complement graph $\overline{G}$ is $(n-k-1)$-regular and $\lambda_{\min}(\overline{G})=-1-\lambda_2(G)\geq\alpha.$ Each connected component of $\overline{G}$ has the same valency $n-k-1\geq d_\theta$, and its least eigenvalue is at least $\alpha$. Hoffman's theorem applies to each component, so every component has least eigenvalue at least $-2$. Therefore, $\lambda_{\min}(\overline{G})\geq-2$, and hence $\lambda_2(G)\leq1$.
\end{proof}

For $1\leq\theta<\sqrt2$, let $k_0(\theta)=\max\{11,d_\theta-1\}$. Suppose that $k\geq k_0(\theta)$, and let $G$ be an extremal graph of order $n=v(k,\theta)$. The crown graph gives $n\geq2k+2$. Hence,
\[
n-k-1\geq k+1\geq d_\theta.
\]
Proposition~\ref{prop:hoffmanjump}, together with $v(k,1)=2k+2$ \cite{Sebisiamvktheta}, gives $v(k,\theta)=2k+2$.

We now extend the conclusion to $1\leq\theta<2$ and determine the graph.
\begin{theorem}\label{thm:jump}
Fix $1\leq\theta<2$. For all sufficiently large $k$, every connected $k$-regular graph with $n\geq2k+2$ and $\lambda_2\leq\theta$ is the crown graph. In particular, its second largest eigenvalue is $1$.
\end{theorem}
\begin{proof}
Suppose that there are connected $k_i$-regular graphs $G_i$ that are not crown graphs, with $k_i\to\infty$, orders $n_i\geq2k_i+2$, and $\lambda_2(G_i)\leq\theta$. Since $1\leq\theta<2$, we have
\[
\left\lfloor1+\frac{\theta^2}{4}\right\rfloor=1.
\]
Theorem~\ref{thm:exactT} therefore gives $n_i=2k_i+2$ for all sufficiently large $i$.

Let $F_i$ be the switched graphs supplied by Lemma~\ref{lem:switch}, and let $d_0$ be a common upper bound on their maximum degrees. This constant depends only on $\theta$. We may assume that $k_i>d_0$. By \eqref{eq:compression},
\[
z^TA_{F_i}z\leq\theta\|z\|^2
\qquad\text{whenever }\langle z,\one\rangle=0.
\]

We first show that $\rho(F_i)\leq\theta$ for all sufficiently large $i$. Let $H$ be a subgraph of $F_i$ on $q$ vertices, and let $x$ be a unit Perron vector of $H$, extended by zero to $F_i$. Let $c=\langle x,\one\rangle,$ $z=x-\frac{c}{n_i}\one.$
Then $c^2\leq q$, $\langle z,\one\rangle=0$, and $\|z\|^2=1-c^2/n_i\leq1$. Since $H$ is a subgraph and $x$ is nonnegative, we have $x^TA_{F_i}x\geq\rho(H)$. Also, $A_{F_i}\one\leq d_0\one$ coordinatewise. Therefore,
\[
\rho(H)-\frac{2d_0q}{n_i}
\leq z^TA_{F_i}z
\leq\theta\|z\|^2
\leq\theta.
\]
In particular,
\[
\rho(H)\leq\theta+\frac{2d_0q}{n_i}.
\]

Choose $\ell$ such that the path $P_\ell$ on $\ell$ vertices satisfies
\[
\rho(P_\ell)=2\cos\frac{\pi}{\ell+1}>\theta.
\]
The preceding estimate shows that $F_i$ cannot contain $P_\ell$ as a subgraph for all sufficiently large $i$. Every vertex in a connected component of $F_i$ is then within distance $\ell-2$ of a fixed vertex of that component. Since the maximum degree is at most $d_0$, each component has order at most
\[
B=1+d_0+\cdots+d_0^{\ell-2}.
\]
There are only finitely many connected graphs of order at most $B$. Applying the same estimate to each of these graphs whose spectral radius exceeds $\theta$ excludes all of them once $i$ is sufficiently large. Hence $\rho(F_i)\leq\theta<2$.

Let $X_i,Y_i$ be the switching cut, and let
\[
t_i=\frac{|X_i|-|Y_i|}{2}.
\]
Since $n_i=2k_i+2$, the number $t_i$ is an integer. Lemma~\ref{lem:cut}, with $s=1$, gives
\[
\overline d(F_i)\geq1+2|t_i|-\frac{2d_0^2}{n_i}.
\]
If $t_i\neq0$, the right-hand side is greater than $2$ for all sufficiently large $i$. This contradicts
\[
\overline d(F_i)\leq\rho(F_i)\leq\theta<2.
\]
Thus $t_i=0$ for all sufficiently large $i$.

Let $w$ have value $1$ on $X_i$ and $-1$ on $Y_i$. The degree formulas in Lemma~\ref{lem:cut}, with $s=1$ and $t_i=0$, give
\[
A_{F_i}w=-w.
\]
Every cycle and $K_{1,4}$ has spectral radius $2$. Since $\rho(F_i)<2$, the graph $F_i$ contains neither, so it is a forest of maximum degree at most $3$. At every vertex $v$, the equation
\[
\sum_{u\sim v}w_u=-w_v
\]
makes $d_{F_i}(v)$ odd. Thus every vertex has degree $1$ or $3$.

A leaf has sign opposite to its neighbor. A vertex of degree $3$ has exactly one neighbor of its own sign, and that neighbor cannot be a leaf. Therefore, any degree-$3$ vertex has a degree-$3$ neighbor. Two adjacent degree-$3$ vertices in a forest contain the double star with two leaves at each center. This double star has spectral radius $2$, since assigning weight $1$ to its centers and $1/2$ to its leaves gives a positive eigenvector with eigenvalue $2$. This again contradicts $\rho(F_i)<2$.

Consequently, $F_i$ is a perfect matching. The equation $A_{F_i}w=-w$ shows that every matching edge crosses $X_i,Y_i$. Both parts have size $k_i+1$, and switching back gives the crown graph. This contradicts the choice of $G_i$.
\end{proof}

\begin{corollary}
For each fixed $1<\theta<2$, there exists $k_0(\theta)$ such that, if $G$ is connected and $k$-regular with $k\geq k_0(\theta)$ and $|V(G)|\geq2k+2$, then
\[
\lambda_2(G)>1\quad\Longrightarrow\quad\lambda_2(G)>\theta.
\]
\end{corollary}

The condition $\theta<2$ is sharp. For $k\geq3$, the bipartite complement of $C_{2(k+2)}$ is connected with valency $k$, order $2k+4$, and has second largest eigenvalue
\[
2\cos\frac{\pi}{k+2}.
\]
These eigenvalues lie in $(1,2)$ and tend to $2$ as $k\to\infty$. Thus, the conclusion of Theorem~\ref{thm:jump} fails at $\theta=2$.

\begin{corollary}\label{cor:least-classification}
Fix $1\leq\theta<2$. There exists an integer $k_1(\theta)$
such that every $k$-regular graph $G$ of order $n\leq 2k$ satisfying $n-k-1\geq k_1(\theta)$ and $\lambda_{\min}(G)\geq-1-\theta$ consists of two disjoint copies of $K_k$ joined by a perfect matching. In particular, $n=2k$ and $\lambda_{\min}(G)=-2$.
\end{corollary}

\begin{proof}
Choose $k_1(\theta)>\theta$ so that Theorem~\ref{thm:jump}
applies to every degree at least $k_1(\theta)$.
Let $H=\overline G$ and $q=n-k-1$. Then $H$ is $q$-regular and
\[
\lambda_2(H)=-1-\lambda_{\min}(G)\leq\theta.
\]
Since $q>\theta$, the graph $H$ is connected.
Also, $n\leq2k$ gives
\[
n\geq2(n-k-1)+2=2q+2.
\]
Theorem~\ref{thm:jump} shows that $H$ is the crown graph on
$2q+2$ vertices. Hence $k=q+1$, and taking its complement gives
two disjoint copies of $K_k$ joined by a perfect matching.
Finally, $\lambda_2(H)=1$, so $\lambda_{\min}(G)=-2$.
\end{proof}

\section{The value of $v(k,\sqrt2)$ for $k\geq13$}\label{sec:sqrt2}
Theorem~\ref{thm:exactT} gives $v(k,\sqrt2)=2k+2$ for all sufficiently large $k$. To prove this for $k\geq13$, we use the classification of Gossett, Jiang, Teets, and Wellner \cite{GJTW}. They show that a graph of order greater than $28$ with the least Seidel eigenvalue at least $-1-2\sqrt2$ is switching equivalent to a disjoint union of isolated vertices, edges, and paths on three vertices. Their classification is computer-assisted. We use the following lemma to determine the regular graphs in these switching classes.

\begin{lemma}\label{lem:cherries}
Let $G$ be a $k$-regular graph with $k\geq2$ and $n\geq2k+2$. If $G$ is switching equivalent to a disjoint union of $P_1,P_2,P_3$, then $G$ is the crown graph.
\end{lemma}

\begin{proof}
Let $F$ be a disjoint union of $P_1,P_2,P_3$ obtained by switching $G$ with respect to a partition $X,Y$. Let $D$ be diagonal with entries $1$ on $X$ and $-1$ on $Y$, and let
\[
w=D\one,\qquad s=\frac n2-k\geq1,\qquad
t=\frac{|X|-|Y|}{2}.
\]
Since $G$ is $k$-regular,
\[
S_G\one=(n-1-2k)\one=(2s-1)\one.
\]
Using $S_F=DS_GD$ and $Dw=\one$, we get $S_Fw=(2s-1)w.$
Also, $Jw=2t\one$, so $S_Fw=2t\one-w-2A_Fw.$
Therefore,
\begin{equation}\label{eq:signswitch}
(A_F+sI)w=t\one.
\end{equation}

For a vertex $v$, \eqref{eq:signswitch} gives
\[
s w_v+\sum_{u\sim_F v}w_u=t.
\]
If $v$ is isolated and $w_v=a\in\{1,-1\}$, then $t=sa$.
On an edge with signs $(a,a)$, we get $t=(s+1)a$.
On an edge with signs $(a,-a)$, the two equations give $t=(s-1)a=-(s-1)a,$ so $s=1$ and $t=0$.

Now consider a component $P_3$. Its leaf equations show that the two leaves have the same sign $a$. The center cannot also have sign $a$, since the leaf and center equations would give $t=(s+1)a$ and $t=(s+2)a$, respectively. Thus, its signs are $(a,-a,a)$, and $t=(s-1)a=(2-s)a.$ Hence, $s=3/2$ and $t=a/2$.

Suppose that $F$ contains a component $P_3$. Then $s=3/2$ and $|t|=1/2$. The equations above exclude isolated vertices and both types of edges. Thus, every component is a $P_3$. Moreover, every component has the same sign pattern $(a,-a,a)$, since $t=a/2$. If there are $r$ components, then $2t=\one^Tw=ra.$ Since $2t=a$, we get $r=1$. This gives $n=3$ and $k=n/2-s=0$, a contradiction.

Suppose next that $F$ contains an edge with opposite signs. Then $s=1$ and $t=0$. An isolated vertex would require $t=\pm1$, and a same-sign edge would require $t=\pm2$. Thus, every component is an opposite-sign edge.

If neither of these cases occurs, the components are isolated vertices or same-sign edges. These two types cannot occur together, since they require $|t|=s$ and $|t|=s+1$, respectively. If all vertices are isolated, their signs agree and switching changes nothing, giving $k=0$. If all components are same-sign edges, their signs again agree throughout the graph, giving $k=1$. Both contradict $k\geq2$.

Consequently, $F$ is a perfect matching with one endpoint of each edge in $X$ and the other in $Y$. Thus $|X|=|Y|=n/2$. Switching back gives the complete bipartite graph between $X$ and $Y$ with this perfect matching deleted. Since $s=1$, we have $n/2=k+1$, so $G$ is the crown graph.
\end{proof}
\begin{theorem}\label{thm:sqrt2}
For every $k\geq13$, $v(k,\sqrt2)=2k+2.$
\end{theorem}
\begin{proof}
The crown graph gives the lower bound. Suppose that a connected $k$-regular graph $G$ has order $n>2k+2$ and $\lambda_2(G)\leq\sqrt2$. On $\one^\perp$, the Seidel eigenvalues are $-1-2\lambda_i(G)\geq-1-2\sqrt2$. On $\one$, the Seidel eigenvalue is $n-1-2k>1$. Thus,
\[\lambda_{\min}(S_G)\geq-1-2\sqrt2.\]
Since $k\geq13$, we have $n\geq29$. The classification in \cite{GJTW} gives that $G$ is switching equivalent to a disjoint union of $P_1, P_2, P_3$. Lemma~\ref{lem:cherries} now shows that $G$ is a crown graph and has order $2k+2$, a contradiction.
\end{proof}

\section{Further questions}\label{sec:questions}
Theorem~\ref{thm:exactT} proves an upper bound for all sufficiently large $k$. We ask for the least degree from which this bound holds.
\begin{question}
For fixed $\theta\geq1$, determine the least integer $K(\theta)$ such that
\[
v(k,\theta)\leq2k+2\left\lfloor1+\frac{\theta^2}{4}\right\rfloor
\qquad(k\geq K(\theta)).
\]
\end{question}

The proof of Corollary~\ref{cor:noncritical} requires $\theta>2\sqrt{L-1}$. This argument does not cover $\theta=2\sqrt{L-1}$, which leads to the following question.
\begin{question}
For each integer $d\geq3$, is $v(k,2\sqrt{d-1})=2k+2d$ for every sufficiently large $k$?
\end{question}
Our construction proves equality along an infinite sequence of degrees. The existence of simple $d$-regular bipartite Ramanujan graphs with every sufficiently large part size would give an affirmative answer. The infinite families we used do not by themselves give all part sizes.

For $1\leq\theta<2$, Theorem~\ref{thm:jump} shows that the crown graph is the unique extremal graph for all sufficiently large $k$. We ask for a similar classification when $\theta\geq2$.
\begin{question}
For fixed $\theta\geq2$, classify the connected $k$-regular graphs of order $2k+2\lfloor1+\theta^2/4\rfloor$ with second largest eigenvalue at most $\theta$, for sufficiently large $k$.
\end{question}

Another question is to determine $v(k,\theta)$ when $\theta$ depends on $k$. For example, Theorem~\ref{thm:quadratic} requires $k>\theta^3$, so it does not apply to $\theta=\sqrt k$. The known bound
\[
v(k,\sqrt k)\leq2k^2+k^{3/2}-k-\sqrt k+1
\]
is discussed in \cite{CGNX}. It would be interesting to improve this bound and determine further values of $v(k,\sqrt k)$.

\section*{Acknowledgement}
The author thanks Igor Balla for a helpful discussion during the ILAS conference at Virginia Tech in May 2026. The author also thanks the organizing committee of ILAS 2026 for financial support through the NSF Early Career Funds.

\end{document}